\documentclass[leqno,12pt]{article}
 \usepackage{amssymb,amsmath}
\newtheorem{theorem}{Theorem}

\newtheorem{proposition}[theorem]{Proposition}
\newtheorem{lemma}[theorem]{Lemma}

\newtheorem{example}[theorem]{Example}

\newcommand{\R}{\mathbb{R}}

\newcommand{\Q}{\mathbb{Q}}

\newcommand{\Sf}{\mathbb{S}}

\newcommand{\Le}{\mathbb{L}}

\newcommand{\Hy}{\mathbb{H}}

\newcommand{\spa}{\mbox{span}}

\newcommand{\rank}{\mbox{rank}}

\newcommand{\po}{{\hspace*{-1ex}}{\bf .  }}

\def\a{\alpha}

\def\<{\langle}
\def\n{\nabla}
\def\>{\rangle}
\def\a{\alpha}

\def\d{\partial}
\def\bea{\begin{eqnarray*} }
\def\eea{\end{eqnarray*} }
\def\be{\begin{equation} }
\def\ee{\end{equation} }

\def\proof{\noindent{\it Proof: }}
\def\qed{\ifhmode\unskip\nobreak\fi\ifmmode\ifinner
\else\hskip5 pt \fi\fi\hbox{\hskip5 pt \vrule width4 pt
height6 pt  depth1.5 pt
\hskip 1pt }}
\begin{document}
\title{ Umbilical surfaces of products of space forms}

\author { Jaime Orjuela\footnote{Supported by CAPES PNPD grant 02885/09-3} and Ruy Tojeiro\footnote{
Partially supported by CNPq grant 	311800/2009-2  and FAPESP grant
2011/21362-2.}}
\date{}
\maketitle
\begin{abstract}
We give a complete  classification of  umbilical surfaces of arbitrary  codimension of a product $\Q_{k_1}^{n_1}\times \Q_{k_2}^{n_2}$ 
of space forms whose curvatures satisfy $k_1+k_2\neq 0$. 
\end{abstract}

\let\thefootnote\relax\footnote{$\!\!\!\!\!\!\!\!\!\!\!$2010 Mathematics Subject Classification. Primary 53B25 Secondary 53C40.}

\section{Introduction}  
A submanifold of a Riemannian manifold is \emph{umbilical} if it is equally curved in all tangent directions. More precisely, an isometric immersion $f\colon\, M^m\to \tilde{M}^n$ between Riemannian manifolds is umbilical if there exists a normal vector field $H$ along $f$ such that its second fundamental form $\alpha_f\in \mbox{Hom}(TM\times TM, N_fM)$
with values in the normal bundle satisfies 
$$\alpha_f(X,Y)=\<X,Y\>H\,\,\,\mbox{for all}\,\,\,X,Y\in \mathfrak{X}(M).$$

The classification of umbilical submanifolds of space forms is very well known. For a general symmetric space $N$, it was shown by Nikolayevsky (see Theorem $1$ of \cite{al}) that any umbilical submanifold of 
 $N$ is an umbilical submanifold of a product of space forms totally geodesically embedded in $N$. This makes the classification of umbilical submanifolds of a product of space forms an important problem.
For submanifolds of dimension $m\geq 3$ of a product $\Q_{k_1}^{n_1}\times \Q_{k_2}^{n_2}$ of space forms  
whose curvatures satisfy $k_1+k_2\neq 0$, the problem was reduced in \cite{mt1} to the classification of $m$-dimensional umbilical submanifolds  with codimension two of $\Sf^n\times \R$ and $\Hy^n\times \R$, where  $\Sf^n$ and $\Hy^n$ stand for the sphere and hyperbolic space, respectively. The case of  $\Sf^n\times \R$ (respectively, $\Hy^n\times \R$) was carried out in \cite{mt2} (respectively, \cite{mt3}), extending previous results in \cite{st} and \cite{vv} (respectively, \cite{ckv}) for  hypersurfaces.

In this paper we extend the results of \cite{mt1} to the surface case.  In this case, the argument in one of the steps of the proof for the higher dimensional case  (see Lemma~$8.2$ of \cite{mt1}) does not apply, and requires more elaborate work. This is carried out in Lemma \ref{le:kerS=0} below, which shows that the difficulty is due to the existence of new interesting families of examples in the surface case. Indeed, our main result (see Theorem \ref{thm:umb} below) states that, in addition to the examples that appear already in higher dimensions, there are precisely two distinct two--parameter families
 of complete embedded flat umbilical surfaces that lie substantially in $\Hy_{k}^{3}\times \R^2$ and $\Hy_{k_1}^{3}\times \Hy_{k_2}^{3}$, respectively. These are discussed in Section~$3$.

\section{Preliminaries}

 Let  $f\colon\,M\to \Q_{k_1}^{n_1}\times \Q_{k_2}^{n_2}$ be  an isometric immersion of a Riemannian manifold.
 We always assume that $M$ is connected.
  Denote by ${\cal R}$  and ${\cal R}^\perp$ the curvature tensors of the tangent and normal bundles
  $TM$ and $N_f M$, respectively,  by $\alpha=\alpha_f\in \Gamma(T^*M\otimes T^*M\otimes N_f M)$ the second fundamental form of $f$ and by $A_\eta=A^f_\eta$ its shape operator  in the normal direction $\eta$, given by  $$\<A_\eta X,Y\>=\<\alpha(X,Y),\eta\>$$ for all $X,Y\in \mathfrak{X}(M)$. Set
 $$L=L^{{f}}:=\pi_2\circ {f}_*\in \Gamma(T^*M\otimes T\Q_{k_2}^{n_2})\,\,\,\mbox{and}\,\,\,K=K^{{f}}:=\pi_2|_{N_fM}\in \Gamma((N_fM)^*\otimes T\Q_{k_2}^{n_2}),$$
 where $\pi_i\colon\, \Q_{k_1}^{n_1}\times \Q_{k_2}^{n_2}\to \Q_{k_i}^{n_i}$ denotes the canonical projection, $1\leq i\leq 2$, and 
by abuse of notation also  its derivative, which we  regard as a  section of $T^*({\Q_{k_1}^{n_1}\times \Q_{k_2}^{n_2}})\otimes T\Q_{k_i}^{n_i}$.

\subsection{The fundamental equations}

The tensors $R\in \Gamma(T^*M\otimes TM)$, $S\in \Gamma(T^*M\otimes N_fM)$ and $T\in \Gamma((N_fM)^*\otimes N_fM)$ given by
\be\label{tensorsa} R=L^tL,\,\,\,\,S=K^tL\,\,\,\,\mbox{and}\,\,\,\,T=K^tK,
\ee
or equivalently, by
$$L=f_*R+S\,\,\,\,\mbox{and}\,\,\,\,K=f_*S^t+T,$$
  were introduced in  \cite{ltv} (see also \cite{mt1}), where
   they were shown to
  satisfy the algebraic relations
\be\label{tensors}
S^tS=R(I-R),\,\,\,\,\,TS=S(I-R)\,\,\,\,
\,\mbox{and}\,\,\,\,\,\,
SS^t=T(I-T),
\ee
as well as the 
differential equations
\be\label{eq:pid1}(\n_XR)Y=A_{SY}X+S^t\a(X,Y),\ee
\be\label{eq:pid2}(\n_XS)Y=T\a(X,Y)-\a(X,RY)\ee
and
\be\label{eq:pid4}(\n_X T)\xi=-SA_\xi X-\a(X,S^t\xi)\ee
for all   $X, Y\in \mathfrak{X}(M)$ and all $\xi \in \Gamma(N_fM)$.
In particular, from the first and third equations of (\ref{tensorsa}) and (\ref{tensors}), respectively, it follows that $R$ and $T$ are  nonnegative operators whose eigenvalues lie in $[0,1]$. 

 The Gauss, Codazzi and Ricci equations of $f$ are, respectively,
 \be\label{eq:gauss}
 {\cal R}(X,Y)Z=(k_1(X\wedge Y-X\wedge R Y-R X\wedge Y) +\kappa R X\wedge R Y)Z
 +A_{\alpha(Y,Z)}X -A_{\alpha(X,Z)}Y,\ee
 \be\label{eq:codazzi}
 (\nabla^\perp_X \alpha)(Y,Z)-(\nabla^\perp_Y \alpha)(X,Z)=\<k_1 X-\kappa RX,Z\>SY-\<k_1 Y-\kappa RY,Z\>SX\ee
 and
 \be\label{eq:ricci}
 {\cal R}^\perp(X,Y)\eta=  \alpha(X,A_\eta Y) -\alpha(A_\eta X,Y)+\kappa(SX\wedge SY)\eta,
 \ee
 where $\kappa=k_1+k_2$.

 \subsection{The flat underlying space}
 
 In order to study isometric immersions $f\colon\,M\to \Q_{k_1}^{n_1}\times \Q_{k_2}^{n_2}$,  it is useful to consider their compositions   $F=h\circ f$ with the 
 canonical isometric embedding
\be\label{eq:j}
h\colon\,\Q_{k_1}^{n_1}\times \Q_{k_2}^{n_2}\to\R_{\sigma(k_1)}^{N_1}\times \R_{\sigma(k_2)}^{N_2}=\R_{\mu}^{N_1+N_2}.
\ee
 Here, for $k\in \R$ we set  $\sigma(k)=1$ if $k<0$ and  $\sigma(k)=0$ otherwise, and as a subscript of an Euclidean space it means the index of the corresponding flat metric.
 Also, we denote  $\mu=\sigma(k_1)+\sigma(k_2)$, $N_i=n_i+1$ if $k_i\neq 0$ and  $N_i=n_i$   otherwise, in which case $\Q_{k_i}^{n_i}$  stands for $\R^{n_i}$.

Let $\tilde\pi_i\colon\,\R_{\mu}^{N_1+N_2}\to \R_{\sigma(k_i)}^{N_i}$, $1\leq i\leq 2$, denote the canonical projection. Then, the normal space of $h$ at each point $z\in \Q_{k_1}^{n_1}\times \Q_{k_2}^{n_2}$ is spanned by $k_1\tilde\pi_1(h(z))$ and $k_2\tilde\pi_2(h(z))$, and its second fundamental form  is given by
\be\label{eq:sffj} \alpha_h(X,Y)=-k_1\<\pi_1 X, Y\>\tilde\pi_1\circ h-k_2\<\pi_2X, Y\>\tilde\pi_2\circ h.\ee
Therefore,  if $k_i\neq 0$,  $1\leq i\leq 2$,   
 then, setting $r_i=|k_i|^{-1/2}$, the unit  vector field  $\nu_i=\nu_i^{F}=\frac{1}{r_i}\tilde\pi_i\circ F$ is  normal  to $F$  and
   we have 
$$\tilde{\nabla}_X\nu_1=\frac{1}{r_1}\tilde\pi_1F_*X=\frac{1}{r_1}(F_*X-h_*LX)=\frac{1}{r_1}(F_*(I-R)X-h_*SX)$$
and
$$\tilde{\nabla}_X\nu_2=\frac{1}{r_2}\tilde\pi_2F_*X=\frac{1}{r_2}h_*LX=\frac{1}{r_2}(F_*RX+h_*SX),$$
where $\tilde \n$ stands for the derivative in $\R_{\mu}^{N_1+N_2}$.
 Hence
 \be\label{eq:nu}^F\nabla_X^\perp \nu_1=-\frac{1}{r_1}h_*SX, \,\,\,\,\,\,\,\,\,\,\,\,\,\,A^{F}_{\nu_1} = -\frac{1}{r_1}(I-R),\ee
 
 \be\label{eq:nu'}^F\nabla_X^\perp \nu_2=\frac{1}{r_2}h_*SX\,\,\,\,\,\,\,\,\mbox{and}\,\,\,\,\,\,\,\,\,A^{F}_{\nu_2} = -\frac{1}{r_2}R.\ee

\subsection{Reduction of codimension}  

 An isometric immersion $f\colon\, M^m\to \Q_{k_1}^{n_1}\times \Q_{k_2}^{n_2}$ is said to \emph{reduce codimension on the left by $\ell$} if there exists a totally geodesic inclusion $j_1\colon\, \Q_{k_1}^{m_1}\to \Q_{k_1}^{n_1}$, with $n_1-m_1=\ell$, and an isometric immersion $\bar f\colon\, M^m\to \Q_{k_1}^{m_1}\times \Q_{k_2}^{n_2}$ such that $f=(j_1\times id)\circ \bar f$. Similarly one defines what it means by $f$ reducing codimension \emph{on the right}.

We will need the  following result  from \cite{mt1} on reduction of codimension. In the statement, $U$ and $V$ stand for $\ker T$ and $\ker (I-T)$, respectively. Notice that the third equation in (\ref{tensors}) implies that   $S(TM)^\perp$ splits orthogonally as $S(TM)^\perp=U\oplus V$, with $U=(I-T)(S(TM)^\perp)$ and $V=T(S(TM)^\perp)$. Also, 
given an isometric immersion $f\colon\, M\to \tilde M$ between Riemannian manifolds, its \emph{first normal space} at $x\in M$ is the subspace $N_1(x)$ of $N_fM(x)$ spanned by the image of its second fundamental form at $x$.

\begin{proposition}\po\label{cor:redcod} Let $f\colon\,M^m\to \Q_{k_1}^{n_1}\times \Q_{k_2}^{n_2}$ be an isometric immersion. If $U\cap N_1^\perp$ (respectively, $V\cap N_1^\perp$ ) is a vector subbundle of $N_fM$ with rank $\ell$  satisfying $\nabla^\perp (U\cap N_1^\perp)\subset N_1^\perp$ (respectively, $\nabla^\perp (V\cap N_1^\perp)\subset N_1^\perp$), then $f$ reduces codimension on the left (respectively, on the right) by $\ell$.
\end{proposition} 

\subsection{Frenet formulae for space-like curves in $\R_{1}^4$}

We briefly recall the definition of the Frenet curvatures and the Frenet frame of a unit-speed space-like curve $\gamma\colon I\to \R_{1}^4$ 
in the four dimensional Lorentz space, as well as the coresponding Frenet formulae, which will be needed in the sequel.

Thus, we assume that $t(s)=\gamma'(s)$ satisfies $\<t(s), t(s)\>=1$ for all $s\in I$. 
Assume first that $\<\gamma''(s), \gamma''(s)\>\neq 0$ for all $s\in I$. Define
$\hat{k}_1(s)=\|\gamma''(s)\|=|\<\gamma''(s), \gamma''(s)\>|^{1/2}$ 
and $n_1(s)=\gamma''(s)/\hat{k}_1(s)$ for all $s\in I$. Denote $\epsilon_1=\<n_1, n_1\>$. Suppose that $v(s)=n_1'(s)+\epsilon_1\hat{k}_1(s)t(s)$ satisfies
$\<v(s), v(s)\>\neq 0$ for all $s\in I$. Define $\hat{k}_2(s)=\|v(s)\|$ and $n_2(s)=v(s)/\hat{k}_2(s)$.  Let $n_3(s)$ be chosen so that 
$\{t(s), n_1(s), n_2(s),n_3(s)\}$ is a positively-oriented orthonormal basis of $\R_{1}^4$ and set $\epsilon_3=\<n_3,n_3\>$. Then the following Frenet formulae hold, where $\hat{k}_3$ is defined by the third equation:
$$\left\{\begin{array}{l}t'=\hat{k}_1n_1,\\
n_1'=-\epsilon_1\hat{k}_1t+\hat{k}_2n_2\\
n_2'=\epsilon_3\hat{k}_2n_1+\hat{k}_3n_3\\
n_3'=\epsilon_1\hat{k}_3n_2.
\end{array}\right.$$

Lesser known are the formulae in the case in which  $\gamma''(s)$ is a nonzero light-like vector everywhere, i.e., $\tilde{n}_1(s)=\gamma''(s)$ satisfies 
$\<\tilde{n}_1(s), \tilde{n}_1(s)\>=0$ for all $s\in I$. We carry them out in more detail below. 

First notice that $\<t, \tilde{n}_1\>=0$. Here, and in the next computations, we drop the ``s"
for simplicity of notation and understand that all equalities hold for all $s\in I$. Thus, 
$$\<\tilde{n}_1', t\>=-\<t',\tilde{n}_1\>=-\<\tilde{n}_1,\tilde{n}_1\>=0.$$
Moreover, $\<\tilde{n}_1', \tilde{n}_1\>=0$, hence $\tilde{n}_1'$ is space-like. Define $\tilde{k}_1=\|\tilde{n}_1'\|$ and $\tilde{n}_2$ by $\tilde{n}_1'=\tilde{k}_1\tilde{n}_2.$
Now let $\tilde{n}_3\in \{t, \tilde{n}_2\}^\perp$ be the unique vector such that
$$\<\tilde{n}_3, \tilde{n}_3\>=0\,\,\,\,\mbox{and}\,\,\,\,\<\tilde{n}_1, \tilde{n}_3\>=1,$$
that is, $\{\tilde{n}_1, \tilde{n}_3\}$ is a pseudo-othonormal basis of the time-like plane $\{t, \tilde{n}_2\}^\perp$.
Since 
$$\<\tilde{n}_2', t\>=-\<\tilde{n}_2, t'\>=-\<\tilde{n}_2, \tilde{n}_1\>=0$$
and 
$$\<\tilde{n}_2', \tilde{n}_1\>=-\<\tilde{n}_2, \tilde{n}_1'\>=-\tilde{k}_1,$$
we have 
$$\tilde{n}_2'=\<\tilde{n}_2', \tilde{n}_1\>\tilde{n}_3+\<\tilde{n}_2', \tilde{n}_3\>\tilde{n}_1=-\tilde{k}_1\tilde{n}_3-\tilde{k}_2\tilde{n}_1,$$
where
$$\tilde{k}_2=\<\tilde{n}_3', \tilde{n}_2\>.$$
Finally, since
$$0=\<\tilde{n}_3',t\>=\<\tilde{n}_3', \tilde{n}_1\>=\<\tilde{n}_3', \tilde{n}_3\>,$$
we have
$$\tilde{n}_3'=\<\tilde{n}_3', \tilde{n}_2\>\tilde{n}_2=\tilde{k}_2\tilde{n}_2.$$

In summary, for a unit-speed space-like curve $\gamma\colon I\to \R_{1}^4$ with light-like curvature vector $\gamma''$, 
one can define two Frenet curvatures $\tilde{k}_1$ and $\tilde{k}_2$ and a pseudo-orthonormal Frenet frame $\{t, \tilde{n}_1, \tilde{n}_2, \tilde{n}_3\}$
with respect to which the Frenet formulae are
$$
\left\{\begin{array}{l} t'=\tilde{n}_1\vspace{1ex}\\
\tilde{n}_1'=\tilde{k}_1\tilde{n}_2\vspace{1ex}\\
\tilde{n}_2'=-\tilde{k}_2\tilde{n}_1-\tilde{k}_1\tilde{n}_3\vspace{1ex}\\
\tilde{n}_3'=\tilde{k}_2\tilde{n}_2.
\end{array}\right.
$$
In both cases, a unit-speed space-like curve $\gamma\colon I\to \R_{1}^4$ is completely determined by its Frenet curvatures, up
to an isometry of $\R_1^4$.

\section{Flat umbilical surfaces in $\Hy_{k}^{3}\times \R^2$ and $\Hy_{k_1}^{3}\times \Hy_{k_2}^{3}$.}

We present below two families  of complete flat properly embedded umbilical surfaces, the first one  in $\Hy_{k}^{3}\times \R^2$ and the second in $\Hy_{k_1}^{3}\times \Hy_{k_2}^{3}$,  each of which depending on two parameters.

\begin{example}\po\label{ex:1} {\em Let $F\colon \R^2\to \R_{1}^6=\R_{1}^4\times \R^2$, where $\R_{1}^4$ has signature $(-, +, +, +)$, be given by
\be\label{eq:F}F(s,t)=\left(a_1\cosh \frac{s}{c}, a_1\sinh \frac{s}{c}, a_2\cos \frac{t}{c}, a_2\sin \frac{t}{c}, b_1\frac{s}{c}, b_2\frac{t}{c}\right),\ee
with 
\be\label{eq:cond}a_1^2-a_2^2=r^2\,\,\,\,\,\,\,\,\mbox{and}\,\,\,\,\,\,\,\,a_1^2+b_1^2=c^2=a_2^2+b_2^2.\ee
Then $F(\R^2)\subset \Hy_{k}^3\times \R^2$, where  $k=-{1}/{r^2}$, by the first relation in (\ref{eq:cond}).
If $\{e_1, \ldots, e_6\}$ is the orthonormal basis of $\R_1^6$ with respect to which $F$ is given by (\ref{eq:F}),
then the subspaces $V_1$ and $V_2$ of $\Le^6$ spanned by $\{e_1, e_2, e_5\}$ and $\{e_3, e_4, e_6\}$
 can be identified
with $\R_1^3$ and $\R^3$, respectively, and 
$$F=\gamma_1\times \gamma_2\colon \R\times \R=\R^2\to V_1\times V_2=\R_1^3\times \R^3=\R_1^6,$$
where $\gamma_1$ and $\gamma_2$ are the helices in $\R_1^3$ and $\R^3$, respectively, parameterized by
$$\gamma_1(s)=\left(a_1\cosh \frac{s}{c}, a_1\sinh \frac{s}{c},b_1\frac{s}{c}\right)$$
and 
$$\gamma_2(t)=\left(a_2\cos \frac{t}{c}, a_2\sin \frac{t}{c},  b_2\frac{t}{c}\right).$$
By the  relations on 
the right  in (\ref{eq:cond}), both $\gamma_1$ and $\gamma_2$ are unit-speed curves, 
hence $F$ is an isometric immersion. Since $F(\R^2)\subset \Hy_{k}^3\times \R^2$, 
there exists an isometric immersion $f\colon \R^2\to \Hy_{k}^3\times \R^2$ such that 
 $F=h\circ f$, where  $h\colon \Hy_{k}^3\times \R^2\to \R_1^6$ denotes  the inclusion.
 It is easily checked that the second fundamental form of $f$ satisfies
$$\alpha_f\left(\frac{\d}{\d s},\frac{\d}{\d t}\right)=0$$
and
$$
\alpha_f\left(\frac{\d}{\d s}, \frac{\d}{\d s}\right)=H(s,t)=\alpha_f\left(\frac{\d}{\d t}, \frac{\d}{\d t}\right),
$$
where
$$h_*H(s,t)=\frac{ka_1a_2}{c^2}\left(a_2\cosh \frac{s}{c}, a_2\sinh \frac{s}{c}, a_1\cos \frac{t}{c}, a_1\sin \frac{t}{c}, 0,0\right).$$
Hence $f$ is   umbilical  with mean curvature vector field $H$.

In view of (\ref{eq:cond}), one can write
$$a_1^2=r^2\frac{(1-\lambda_1)}{\lambda_2-\lambda_1},\,\,\,\,a_2^2=r^2\frac{(1-\lambda_2)}{\lambda_2-\lambda_1},\,\,\,\,b_1^2=r^2\frac{\lambda_1}{\lambda_2-\lambda_1},\,\,\,\,b_2^2=r^2\frac{\lambda_2}{\lambda_2-\lambda_1},\,\,\,\,c^2=\frac{r^2}{\lambda_2-\lambda_1},$$
with  $0<\lambda_1< \lambda_2<1$.
Then, one can check that the  curvature vector $\gamma_i''$ of  $\gamma_i$, $1\leq i\leq 2$, satisfies
\be\label{eq:curvvector}\<\gamma_i'', \gamma_i''\>=k(\lambda_j-\lambda_i)(1-\lambda_i), \,\,\,1\leq i\neq j\leq 2,\ee
and that  the second Frenet curvature (torsion) of $\gamma_i$ satisfies
\be\label{torsion}\tau_i^2=-k\lambda_i|\lambda_j-\lambda_i|,\,\,\,\,1\leq i\neq j\leq 2.\ee  }
\end{example}

\begin{example}\po \label{ex:2}{\em Let $\R_2^8=\R_1^4\times \R_1^4$ denote Euclidean space of dimension $8$ endowed with 
an inner product of signature $(-, +,+,+,-,+,+, +)$, and let $F\colon \R^2\to \R_2^8$ be given by
\be\label{eq:Fb}F(s,t)=\left(a_1\cosh \frac{s}{c}, a_1\sinh \frac{s}{c}, a_2\cos \frac{t}{c}, a_2\sin \frac{t}{c}, a_3\cosh \frac{t}{d}, a_3\sinh \frac{t}{d}, a_4\cos \frac{s}{d}, a_4\sin \frac{s}{d} \right),\ee
with
\be\label{eq:condb}a_1^2-a_2^2=r_1^2,\,\,\,\,\,\,\,\,a_3^2-a_4^2=r_2^2\,\,\,\,\,\,\,\,\,\mbox{and}\,\,\,\,\,\,\,\,
\frac{a_1^2}{c^2}+\frac{a_4^2}{d^2}=1=\frac{a_2^2}{c^2}+\frac{a_3^2}{d^2}.\ee
The first pair of relations in (\ref{eq:condb}) implies that $F(\R^2)\subset \Hy_{k_1}^3\times \Hy_{k_2}^3\subset \R_1^4\times \R_1^4$, with $k_i=-1/r_i^2$
 for  $1\leq i\leq 2$. If $\{e_1, \ldots, e_4, f_1, \ldots, f_4\}$ is the orthonormal basis of $\R_2^8$ with respect to which $F$ is given by (\ref{eq:Fb}),
then the subspaces $V_1$ and $V_2$ of $\R_2^8$ spanned by $\{e_1, e_2, f_3, f_4\}$ and $\{f_1, f_2, e_3, e_4\}$
 can also be identified
with $\R_1^4$, and 
$$F=\gamma_1\times \gamma_2\colon \R\times \R=\R^2\to V_1\times V_2,$$
where $\gamma_1$ and $\gamma_2$ are the curves  parameterized by
$$\gamma_1(s)= \left(a_1\cosh \frac{s}{c}, a_1\sinh \frac{s}{c},  a_4\cos \frac{s}{d}, a_4\sin \frac{s}{d} \right)$$
and
$$\gamma_2(t)=\left(a_3\cosh \frac{t}{d}, a_3\sinh \frac{t}{d}, a_2\cos \frac{t}{c}, a_2\sin \frac{t}{c}\right).$$
 In view of the second pair of relations in (\ref{eq:condb}), both $\gamma_1$ and $\gamma_2$ are unit-speed curves, 
hence $F$ is an isometric immersion. Since $F(\R^2)\subset \Hy_{k_1}^3\times \Hy_{k_2}^3$, 
there exists an isometric immersion $f\colon \R^2\to \Hy_{k_1}^3\times \Hy_{k_2}^3$ such that 
 $F=h\circ f$, where  $h\colon \Hy_{k_1}^3\times \Hy_{k_2}^3\to \R_2^8$ denotes  the inclusion.
 One can  easily check that the second fundamental form of $f$ satisfies
 $$\alpha_f\left(\frac{\d}{\d s},\frac{\d}{\d t}\right)=0$$
and
$$
\alpha_f\left(\frac{\d}{\d s}, \frac{\d}{\d s}\right)= H(s,t)
=\alpha_f\left(\frac{\d}{\d t}, \frac{\d}{\d t}\right)
$$
where
$$\begin{array}{l} h_*H(s,t)=\frac{k_1a_1a_2}{c^2}\left(a_2\cosh \frac{s}{c}, a_2\sinh \frac{s}{c}, a_1\cos \frac{t}{c}, a_1\sin \frac{t}{c},0,0,0,0\right)+ \vspace{1ex}\\
\frac{k_2a_3a_4}{d^2}\left(0,0,0,0,a_4\cosh \frac{t}{d}, a_4\sinh \frac{t}{d}, a_3\cos \frac{s}{d}, a_3\sin \frac{s}{d}\right).
\end{array}$$
It follows that $f$ is  umbilical with mean curvature vector field $H$.

By the conditions in (\ref{eq:condb}), one can write
$$a_1^2=r_1^2\frac{(1-\lambda_1)}{\lambda_2-\lambda_1},\,\,\,\,a_2^2=r_1^2\frac{(1-\lambda_2)}{\lambda_2-\lambda_1},\,\,\,\,\,
a_3^2=r_2^2\frac{\lambda_2}{\lambda_2-\lambda_1},\,\,\,\,a_4^2=r_2^2\frac{\lambda_1}{\lambda_2-\lambda_1},$$
$$c^2=\frac{r_1^2}{\lambda_2-\lambda_1}\,\,\,\mbox{and}\,\,\,\,\,d^2=\frac{r_2^2}{\lambda_2-\lambda_1},$$
with  $0<\lambda_1< \lambda_2<1$.  Then, the  curvature vector $\gamma_i''$ of the curve  $\gamma_i$, $1\leq i\leq 2$, satisfies
\be\label{eq:curvvectorb}\<\gamma_i'', \gamma_i''\>=(\lambda_i-\lambda_j)(\kappa\lambda_i-k_1), \,\,\,1\leq i\neq j\leq 2,\,\,\,\kappa=k_1+k_2.\ee
If $\kappa\lambda_i-k_1\neq 0$, one can check that $\gamma_i$, $1\leq i\leq  2$,  has constant Frenet curvatures
$\hat{k}^i_\ell$, $1\leq \ell\leq 3$, given by
\be\label{eq:hatk1}(\hat{k}^i_1)^2=|(\lambda_i-\lambda_j)(\kappa\lambda_i-k_1)|,\ee
\be\label{eq:hatk2}(\hat{k}^i_2)^2=\frac{\kappa^2|\lambda_i-\lambda_j|\lambda_i(1-\lambda_i)}{|\kappa\lambda_i-k_1|}\ee
and
\be\label{eq:hatk3}(\hat{k}^i_3)^2=\frac{k_1k_2|\lambda_i-\lambda_j|}{|\kappa\lambda_i-k_1|},\,\,\,1\leq j\neq i\leq 2.\ee

If $\kappa\lambda_i-k_1= 0$, that is, the curvature vector of $\gamma_i$ is light-like, then one can check that $\gamma_i$ has constant Frenet curvatures $\tilde{k}^i_1$ and $\tilde{k}^i_2$, $1\leq i\leq 2$ (see Subsection $2.4$), given by
\be\label{eq:tildek1}(\tilde{k}^i_1)^2=\frac{k_1k_2(\kappa\lambda_j-k_1)^2}{\kappa^2},\,\,\,1\leq j\neq i\leq 2,\ee
and
\be\label{eq:tildek2}(\tilde{k}^i_2)^2=\frac{(k_1-k_2)^2}{4k_1k_2},\,\,\,1\leq i\leq 2.\ee}
\end{example}

It is also easily checked that the isometric immersions in both of the preceding examples have the frame of
coordinate vector fields $\{\frac{\d}{\d s},\frac{\d}{\d t}\}$ as a frame of principal directions for the associated
tensor $R$, with corresponding eigenvalues $\lambda_1$ and $\lambda_2$, respectively.  Moreover, they are clearly injective
and proper, hence embeddings. Therefore, all surfaces in both families are properly embedded and isometric to the plane.

\section{The main step}

 Umbilical submanifolds of $\Q_{k_1}^{n_1}\times \Q_{k_2}^{n_2}$ were studied in \cite{mt1} according to the possible structures of the tensor $S$. When $\ker S=\{0\}$, it
was shown that $R$ must be a constant multiple of the identity tensor whenever the dimension of the submanifold is at least three (see \cite{mt1}, Lemma $8.2$), which corresponds to case $(i)$ in the statement of Lemma \ref{le:kerS=0}  below. We  now show that in the surface case 
the only exceptions are the surfaces of the two families in the preceding section. 

\begin{lemma}\po\label{le:kerS=0} Let $f\colon\,M^2\to \Q_{k_1}^{n_1}\times \Q_{k_2}^{n_2}$, $k_1+k_2\neq 0$, be an umbilical isometric immersion. Assume that   $\ker S=\{0\}$ at some point $x\in M^2$. Then  one of the following  holds:
\begin{itemize}
\item[(i)] there exist umbilical isometric immersions $f_i\colon\, M^2\to \Q_{\tilde k_i}^{n_i}$, $1\leq i\leq 2$, with $\tilde{k}_1=k_1\cos^2\theta$ and $\tilde k_2=k_2\sin^2\theta$ for some $\theta\in (0, \pi/2)$,  such that $f=(\cos \theta f_1,\sin\theta f_2)$. 
\item[(ii)] after interchanging the factors, if necessary, we have $k_2=0$, $n_1\geq 3$, $n_2\geq 2$ and $f=j\circ \tilde f$, where $j\colon \Q_{k_1}^{3}\times \R^2\to \Q_{k_1}^{n_1}\times \R^{n_2}$ and $\tilde f\colon M^2\to \Q_{k_1}^{3}\times \R^2$ are isometric immersions such that $j$ is totally geodesic and $\tilde f(M^2)$ is an open subset of a surface as in Example \ref{ex:1}.
\item[(iii)]  $k_i<0$ and  $n_i\geq 3$, $1\leq i\leq 2$,  and $f=j\circ \tilde f$, where 
$j\colon \Q_{k_1}^{3}\times \Q_{k_2}^{3}\to \Q_{k_1}^{n_1}\times \Q_{k_2}^{n_2}$ and  $\tilde f\colon M^2\to \Q_{k_1}^{3}\times \Q_{k_2}^{3}$ are isometric immersions such that $j$ is totally geodesic and $\tilde f(M^2)$ is an open subset of a surface as in Example \ref{ex:2}.
\end{itemize}
\end{lemma}
 
\proof Let ${\cal U}\subset M$ be the maximal  connected open  neighborhood of  $x$ where   $\ker S=\{ 0\}$. 
Let $\lambda_1$ and $\lambda_2$  be the eigenvalues of $R$ on ${\cal U}$. Fix an orthonormal frame $\{X_1,X_2\}$ of eigenvectors of $R$, with $X_i$ associated to $\lambda_i$,  and define $\xi_i:=SX_i$ for $i=1$, $2$. Thus, from (\ref{tensors}) we have
\begin{equation}\label{nnf}
\langle\xi_i,\xi_j\rangle = \langle S^tSX_i,X_j\rangle = \delta_{ij}\lambda_i(1-\lambda_i)
\end{equation}
and
\begin{equation}\label{evT}
T\xi_i=TSX_i=(1-\lambda_i)\xi_i
\end{equation}
for $i$, $j=1,2$. We can write equations (\ref{eq:gauss})--(\ref{eq:ricci})  in the frames $\{X_1, X_2\}$ and $\{\xi_1, \xi_2\}$, in terms of the Gaussian curvature $K$ of $M^2$ and the mean curvature vector $H$ of $f$, as 
\be\label{eq:gauss3}
K=k_1(1-\lambda_1)(1-\lambda_2)+k_2\lambda_1\lambda_2+|H|^2,
\ee
\be\label{eq:codazzi3}
\nabla_{X_i}^{\perp}H=(\kappa \lambda_j-k_1)\xi_i, \,\,\,1\leq i\neq j\leq 2,
\end{equation}
and 
\be\label{eq:ricci3}
\mathcal{R}^{\perp}(X_1,X_2)=\kappa (\xi_1\wedge\xi_2),
\ee
whereas equations  (\ref{eq:pid1})--(\ref{eq:pid4}) become
\begin{equation}
\label{cdt1b}
\left(\nabla_{X_i} R\right)X_j=\langle \xi_j,H\rangle X_i+\delta_{ij} S^tH,
\ee
\be
\label{cdt2b}
\left(\nabla_{X_i} S\right)X_j=\delta_{ij}(T-\lambda_jI)H
\ee
and
\be\label{cdt3b}
\left(\nabla_{X_i} T\right)\xi=-\langle \xi,H\rangle \xi_i-\langle \xi,\xi_i\rangle H
\end{equation}
for $i$, $j=1,2$ and all $\xi\in \Gamma(N_fM^2)$.
Define the Christoffel symbols $\Gamma_{11}^2$ and $\Gamma_{22}^1$ by 
\begin{equation}\label{lcc}
\nabla_{X_1}X_1=\Gamma_{11}^2X_2\,\,\,\,\,\mbox{and}\,\,\,\,\,\,
\nabla_{X_2}X_2=\Gamma_{22}^1X_1.
\end{equation}
Substituting
\begin{align*}
\left(\nabla_{X_i}R\right)X_j &= \nabla_{X_i}RX_j-R\nabla_{X_i}X_j\\
                              &= X_i(\lambda_j)X_j + (\lambda_jI-R)\nabla_{X_i}X_j
\end{align*}
into  (\ref{cdt1b}) yields
\be \label{cdlamda}X_i(\lambda_j)=\delta_{ij}2\langle\xi_i,H\rangle\ee
and
\begin{equation}\label{relation}
\langle\xi_i,H\rangle=(\lambda_j-\lambda_i)\Gamma_{jj}^i,\,\,\,\,1\leq i\neq j\leq 2.
\end{equation}
On the other hand, from    (\ref{cdt2b})  we get
\begin{equation}\label{pnc-f}
\nabla^{\perp}_{X_i}\xi_j=-\Gamma_{ii}^j\xi_i, \,\,\,1\leq i\neq j\leq 2.
\end{equation}
Using  (\ref{eq:codazzi3}),    (\ref{cdlamda}) and (\ref{pnc-f}) we obtain
\begin{align*}
\mathcal{R}^{\perp}(X_1,X_2)H= &\nabla^{\perp}_{X_1}\nabla^{\perp}_{X_2}H-\nabla^{\perp}_{X_2}\nabla^{\perp}_{X_1}H-\nabla^{\perp}_{[X_1,X_2]}H\vspace{1.5ex}\\
                             =& \nabla^{\perp}_{X_1}(\kappa \lambda_1-k_1)\xi_2-\nabla^{\perp}_{X_2}(\kappa \lambda_2-k_1)\xi_1
                             +(\kappa \lambda_2-k_1)\Gamma_{11}^2\xi_1\vspace{1.5ex}\\&-(\kappa \lambda_1-k_1)\Gamma_{22}^1\xi_2\vspace{1.5ex}\\
                             =& \kappa  X_1(\lambda_1)\xi_2+(\kappa \lambda_1-k_1)\nabla^{\perp}_{X_1}\xi_2
                             -\kappa  X_2(\lambda_2)\xi_1-(\kappa \lambda_2-k_1)\nabla^{\perp}_{X_2}\xi_1\vspace{1.5ex}\\
                             &+(\kappa \lambda_2-k_1)\Gamma_{11}^2\xi_1-(\kappa \lambda_1-k_1)\Gamma_{22}^1\xi_2\vspace{1.5ex}\\
                             =& 2\kappa \langle\xi_1,H\rangle\xi_2-(\kappa \lambda_1-k_1)\Gamma_{11}^2\xi_1
                             -2\kappa \langle\xi_2,H\rangle\xi_1+(\kappa \lambda_2-k_1)\Gamma_{22}^1\xi_2\vspace{1.5ex}\\
                             &+(\kappa \lambda_2-k_1)\Gamma_{11}^2\xi_1-(\kappa \lambda_1-k_1)\Gamma_{22}^1\xi_2\vspace{1.5ex}\\
                             =&-\kappa (2\langle\xi_2,H\rangle+(\lambda_1-\lambda_2)\Gamma_{11}^2)\xi_1
                             +\kappa (2\langle\xi_1,H\rangle+(\lambda_2-\lambda_1)\Gamma_{22}^1)\xi_2.
\end{align*}
In view of (\ref{relation}), the above equation becomes
\begin{equation}\label{ricci-H-def}
\mathcal{R}^{\perp}(X_1,X_2)H=-3\kappa (\langle\xi_2,H\rangle\xi_1-\langle\xi_1,H\rangle\xi_2).
\end{equation}
Comparing the preceding equation with 
\begin{equation}\label{ricci-H}
\mathcal{R}^{\perp}(X_1,X_2)H=\kappa (\langle\xi_2,H\rangle\xi_1-\langle\xi_1,H\rangle\xi_2),
\end{equation}
which follows from  (\ref{eq:ricci3}), and using that $\kappa\neq 0$, we get $\langle\xi_1,H\rangle=0=\langle\xi_2,H\rangle$, 
i.e.,
\begin{equation}\label{eq:hperp} H\in \Gamma(S(TM)^{\perp}).
\end{equation}
 In particular, we obtain from   (\ref{cdlamda}) that $\lambda_1$ and $\lambda_2$ assume constant values in $(0,1)$ everywhere on ${\cal U}$. 
If ${\cal U}$ were a proper subset of $M^2$, then $\lambda_1$ and $\lambda_2$ would assume  the same values on the boundary of ${\cal U}$, hence $\lambda_i(1-\lambda_i)\neq 0$ on an open connected neighborhood
of $\bar{{\cal U}}$, $1\leq i\leq 2$, contradicting the maximality of ${\cal U}$ as an open  connected neighborhood of $x$ where   $\ker S=\{ 0\}$. It follows that   ${\cal U}=M^2$. 
Moreover, if $\lambda_1=\lambda_2$ then $f$ is as in $(i)$ by Proposition $5.2$ of \cite{mt1}. From now on we assume that $\lambda_1\neq \lambda_2$. 

Under this assumption,  we obtain from (\ref{relation}) and (\ref{eq:hperp}) that  $\Gamma_{11}^2=0=\Gamma_{22}^1$. In particular, we have $K=0$ everywhere,
and then  (\ref{eq:gauss3}) gives
\begin{equation}\label{eq:H}
|H|^2=-k_1(1-\lambda_1)(1-\lambda_2)-k_2\lambda_1\lambda_2.
\end{equation}

Set $\xi=H$ in (\ref{cdt3b}). By using (\ref{eq:codazzi3}) and (\ref{eq:H}),  we obtain
\begin{equation}\label{ndc-TH}
\nabla^{\perp}_{X_i}TH=k_2\lambda_j\xi_i, \,\,\,\, 1\leq i\neq j\leq 2,
\end{equation}
and similarly
\begin{equation}\label{ndc-(I-T)H}
\nabla^{\perp}_{X_i}(I-T)H=-(1-\lambda_j)k_1\xi_i, \,\,\, 1\leq i\neq j\leq 2.\\
\end{equation}
In particular, bearing in mind (\ref{eq:hperp}) and the fact that $T$ leaves $S(TM)$ invariant, as follows from 
the second equation in (\ref{tensors}), we obtain that  both $TH$ and $(I-T)H$ have constant length on $M^2$, hence either $TH=0$, $TH=H$ or both $TH$ and $(I-T)H$ are nonzero  everywhere.  
Therefore $L_1=U\cap \{H\}^{\perp}=U\cap N_1^\perp$ and $L_2=V\cap \{H\}^{\perp}=V\cap N_1^\perp$ have constant dimensions on $M^2$, which are, accordingly,
 $(\rank\, U -1, \rank\, V)$, $(\rank\, U, \rank\, V\,-1)$ or $(\rank\, U-1, \rank\, V-1)$. 
Moreover,  equations (\ref{eq:codazzi3}) and (\ref{eq:hperp})
imply that $\nabla^{\perp}_{TM}L_{i}\subset \{H\}^{\perp}$ for $i=1,2$. Hence, the assumptions of Proposition \ref{cor:redcod} are satisfied, and we conclude that 
there are three corresponding possibilities for the pairs  $(n_1,n_2)$ of  \emph{substantial} values of $n_1$ and $n_2$: $(3, 2)$, $(2, 3)$ and $(3,3)$.

We first consider the case  $(n_1,n_2)=(3,2)$. This is the case in which $TH=0$, and hence $k_2=0$ by (\ref{ndc-TH}). Thus we have $k_1<0$ from (\ref{eq:H}), 
and we may assume that $f$ takes values in $\Hy_{k}^3\times \R^2$, with $k=k_1<0$.

 Set $F=h\circ f$, where $h\colon \Hy_{k}^3\times \R^2\to \R_1^6$ denotes the inclusion. 
The second fundamental form of $F$ is given by
$$\alpha_F(X, Y)=\<X, Y\>h_*H+\frac{1}{r}\<(I-R)X, Y\>\nu,$$
where $r=(-k)^{-1/2}$ and $\nu=\frac{1}{r}\tilde \pi_1\circ F$. Therefore
\be\label{eq:zi}\alpha_F(X_i, X_j)=\delta_{ij}(h_*H+\frac{1}{r}(1-\lambda_i)\nu):=\delta_{ij}Z_i=\tilde \nabla_{X_j}F_*X_i,\,\,\,1\leq i, j\leq 2.\ee
Notice that
$$\<Z_1, Z_2\>=|H|^2+k(1-\lambda_1)(1-\lambda_2)=0$$
 by (\ref{eq:H}), 
 and that
\be\label{eq:zinorm}\<Z_i, Z_i\>=k(\lambda_j-\lambda_i)(1-\lambda_i),\,\,\,\,\,\,1\leq i\neq j\leq 2.\ee
Moreover, since
$$\tilde{\pi}_2(h_*H)=h_*\pi_2 H=h_*(f_*S^tH+TH)=0,$$
it follows that
$$\tilde{\pi}_2Z_i=0, \,\,\,1\leq i\leq 2.$$
Using (\ref{eq:codazzi3}), we have
\begin{eqnarray*} \tilde \nabla_{X_i}h_*H
&=& h_*\hat \nabla_{X_i}H+\alpha_h(f_*X_i, H)
\vspace{1ex}\\
&=&-F_*A^f_HX_i+h_*\nabla^\perp_{X_i}H+\frac{1}{r}\<\pi_1f_*X_i, H\>\nu\vspace{1ex}\\
&=& -|H|^2F_*X_i-k(1-\lambda_j)h_*\xi_i+\frac{1}{r}\<f_*(I-R)X_i-SX_i,H\>\nu\vspace{1ex}\\
&=&k(1-\lambda_j)\left((1-\lambda_i)F_*X_i-h_*\xi_i\right),\,\,\,\,\,\,1\leq i\neq j\leq 2.
\end{eqnarray*}
On the other hand,
$$ \tilde \nabla_{X_i}\nu
= \frac{1}{r}(F_*(I-R)X_i-h_*SX_i)
= \frac{1}{r}((1-\lambda_i)F_*X_i-h_*\xi_i).
$$
Therefore
\be \label{eq:xizj}\tilde\nabla_{X_i}Z_j=0,\,\,\,\mbox{if}\,\,i\neq j,\ee
and
\be\label{eq:h1}\tilde\nabla_{X_i}Z_i=k(\lambda_i-\lambda_j)((1-\lambda_i)F_*X_i-h_*\xi_i),\,\,\,1\leq i\neq j\leq 2.\ee
Also, using that
$$\nabla_{X_i}^\perp \xi_i=-\frac{1}{|H|^2}\<\nabla_{X_i}H, \xi_i\>H=-\lambda_iH,$$
as follows from (\ref{nnf}), (\ref{eq:codazzi3}) and (\ref{eq:H}), we obtain that 
\begin{eqnarray} \label{eq:h2}\tilde \nabla_{X_i}h_*\xi_j&=&h_*\hat \nabla_{X_i}\xi_j+\alpha_h(f_*X_i, \xi_j)\nonumber\vspace{1ex}\\
&=&-F_*A^f_{\xi_j}X_i+h_*\nabla^\perp_{X_i}\xi_j+\frac{1}{r}\<\pi_1f_*X_i, \xi_j\>\nu\nonumber\vspace{1ex}\\
&=&-\delta_{ij}\lambda_iZ_i, \,\,\,\,1\leq i, j\leq 2.
\end{eqnarray}
It follows from (\ref{eq:zi}), (\ref{eq:xizj}), (\ref{eq:h1}) and (\ref{eq:h2}) that the subspaces $V_i=\spa \{F_*X_i, Z_i, h_*\xi_i\}$, $1\leq i\leq 2$,  are constant. Moreover, they are orthogonal to each other,
hence $\R_1^6$ splits orthogonally as $\R_1^6=V_1\oplus V_2$. 

Since $\Gamma_{11}^2=\Gamma_{22}^1=0$, for each $x\in M^2$ there exists an isometry  $\psi\colon\, W=I_1\times I_2\to U_x$
of a product of open intervals $I_j\subset \R$, $1\leq j\leq 2$, onto an open neighborhood of $x$,  such that $\psi_*\frac{\d}{\d s}=X_1$ 
and $\psi_*\frac{\d}{\d t}=X_2$, where $s$ and $t$ are the standard coordinates on $I_1$ and $I_2$, respectively.  Write $g=F\circ \psi$. In terms of the coordinates $(s,t)$, the fact that $\alpha_F(X_1, X_2)=0$
translates into 
$$\frac{\d^2 g}{\d s\d t}=0,$$
which implies that there exist smooth curves $\gamma_1\colon I_1\to V_1$ and $\gamma_2\colon I_2\to V_2$ such that $g=\gamma_1\times \gamma_2$.
By (\ref{eq:zi}),  (\ref{eq:h1}) and (\ref{eq:h2}), each $\gamma_i$ is a helix in $V_i$  with curvature vector $\gamma_i''=Z_i$ and binormal vector $h_*(\xi_i/|\xi_i|)$, $1\leq i\leq 2$.  It follows from (\ref{eq:zinorm}) that  (\ref{eq:curvvector}) holds for $\gamma_i$, whereas (\ref{eq:h2}) implies that the second Frenet curvature of $\gamma_i$ satisfies
$$\tau_i^2=\frac{\lambda_i^2|\<Z_i, Z_i\>|}{\<\xi_i, \xi_i\>}=-k\lambda_i|\lambda_j-\lambda_i|,\,\,\,\,1\leq i\neq j\leq 2.$$
Therefore, the helices $\gamma_1$ and $\gamma_2$ are precisely, up to congruence, those given in Example~\ref{ex:1}. Moreover, since the curvature vector $Z_i$ along $\gamma_i$  spans a two-dimensional subspace of $V_i$ orthogonal to the axis of $\gamma_i$  and $\tilde \pi_2 Z_i=Z_i$, $1\leq i\leq 2$, it follows that the subspace $\R^2$ in the orthogonal 
decomposition $\R_1^6=\R_1^4\oplus \R^2$ adapted to $\Hy_k^3\times \R^2$ is spanned by the axes of $\gamma_1$ and $\gamma_2$. We conclude that  $g$ is (the restriction to $W$ of ) an isometric immersion as in Example \ref{ex:1}.

We have shown that, for each $x\in M^2$, there exists an open neighborhood $U_x$ of $x$ such that $f(U_x)$ is contained in a 
 surface as in Example~\ref{ex:1} in a totally geodesic $\Q_{k_1}^{3}\times \R^2\subset \Q_{k_1}^{n_1}\times \R^{n_2}$. A standard  
 connectedness argument now shows that $f$ is as in the statement.

The case $(n_1, n_2)=(2,3)$  is entirely similar and leads to the same conclusion with the factors interchanged. Let us consider the case $(n_1, n_2)=(3,3)$,
so we may now assume that $f$ takes values in $\mathbb{Q}^3_{k_1}\times\mathbb{Q}^3_{k_2}$. Here both $TH$ and $(I-T)H$ are nonzero everywhere, 
and  we can choose unit vector fields $\xi_3\in\ker T$,   $\xi_4\in\ker(I-T)$ and write $H=\rho_3\xi_3+\rho_4\xi_4$,
where $\rho_k=\langle\xi_k,H\rangle\not=0$ for $k=3$, $4$. Applying (\ref{cdt3b}) to $\xi=\xi_k$, with $k=3$, $4$,   we get 
 \begin{equation*}
T\nabla^{\perp}_{X_i}\xi_3=\rho_3\xi_i\,\,\,\,\,\,\mbox{and}\,\,\,\,\,\,
(I-T)\nabla^{\perp}_{X_i}\xi_4=-\rho_4\xi_i
\end{equation*}
for $i=1$, $2$. Therefore, for $i=1$, $2$ we obtain
\begin{equation}\label{ncd-imS-perp}
\nabla^{\perp}_{X_i}\xi_3=\frac{\rho_3}{1-\lambda_i}\xi_i\,\,\,\,\,\,\mbox{and}\,\,\,\,\,\,
\nabla^{\perp}_{X_i}\xi_4=-\frac{\rho_4}{\lambda_i}\xi_i.
\end{equation}
Using (\ref{nnf}), the preceding equations yield
\be\label{nxixii}
\nabla^{\perp}_{X_i}\xi_i=-\lambda_i\rho_3\xi_3+(1-\lambda_i)\rho_4\xi_4,\,\,\,\,1\leq i\leq 2.
\ee
On the other hand, we have
\begin{equation}\label{TH-(I-T)H}
(I-T)H=\rho_3\xi_3\,\,\,\,\,\,\mbox{and}\,\,\,\,\,\,
TH=\rho_4\xi_4.
\end{equation}
Thus, combining (\ref{ndc-TH}), (\ref{ndc-(I-T)H}) and (\ref{TH-(I-T)H}) we get
\begin{equation}\label{rho}
\rho_3^2=-k_1(1-\lambda_1)(1-\lambda_2)\,\,\,\,\,\,\mbox{and}\,\,\,\,\,\,
\rho_4^2=-k_2\lambda_1\lambda_2.
\end{equation}
In particular, we must have $k_1, k_2<0$, so  $f$ takes values in $\Hy_{k_1}^3\times \Hy_{k_2}^3$.

Set $F=h\circ f$, where $h\colon \Hy_{k_1}^3\times \Hy_{k_2}^3\to \R_1^4\times \R_1^4=\R_2^8$ denotes the inclusion. 
The second fundamental form of $F$ is given by
$$\alpha_F(X, Y)=\<X, Y\>h_*H+\frac{1}{r_1}\<(I-R)X, Y\>\nu_1+\frac{1}{r_2}\<RX, Y\>\nu_2,$$
where $r_i=(-k_i)^{-1/2}$ and $\nu_i=\frac{1}{r_i}\tilde \pi_i\circ F$,  $1\leq i\leq 2$. Therefore
\be\label{eq:zib}
\alpha_F(X_i, X_j)=\delta_{ij}(h_*H+\frac{1}{r_1}(1-\lambda_i)\nu_1+\frac{1}{r_2}\lambda_i\nu_2):=\delta_{ij}Z_i=\tilde{\nabla}_{X_j}F_*X_i,\,\,\,1\leq i\leq 2.\ee
Notice that
$$\<Z_i, Z_j\>=|H|^2+k_1(1-\lambda_i)(1-\lambda_j)+k_2\lambda_i\lambda_j,\,\,\,\,\,\,1\leq i,  j\leq 2.$$
It follows from (\ref{eq:H}) that 
$$\<Z_1, Z_2\>=0$$
and
\be\label{eq:zinormb}
\<Z_i, Z_i\>=(\lambda_i-\lambda_j)(\kappa\lambda_i-k_1),\,\,\,\,\,\,1\leq i\neq j\leq 2.
\ee
Using (\ref{eq:codazzi3}), we obtain
\begin{eqnarray*} \tilde \nabla_{X_i}h_*H
&=& h_*\hat \nabla_{X_i}H+\alpha_h(f_*X_i, H)
\vspace{1ex}\\
&=&-F_*A^f_HX_i+h_*\nabla^\perp_{X_i}H\vspace{1ex}\\
&=& -|H|^2F_*X_i+(\kappa \lambda_j-k_1)h_*\xi_i.
\end{eqnarray*}
On the other hand,
$$
 \tilde \nabla_{X_i}\nu_1
= \frac{1}{r_1}(F_*(I-R)X_i-h_*SX_i)= \frac{1}{r_1}((1-\lambda_i)F_*X_i-h_*\xi_i)
$$
and 
$$
 \tilde \nabla_{X_i}\nu_2
= \frac{1}{r_2}(F_*RX_i+h_*SX_i)
= \frac{1}{r_2}(\lambda_iF_*X_i+h_*\xi_i).
$$
Using (\ref{eq:H}), it follows that
$$\tilde\nabla_{X_i}Z_j=0,\,\,\,\mbox{if}\,\,i\neq j,$$
whereas
\be\label{eq:h1b}\tilde\nabla_{X_i}Z_i=-\<Z_i,Z_i\>F_*X_i+\kappa (\lambda_j-\lambda_i)h_*\xi_i,\,\,\,1\leq i\neq j\leq 2.\ee
Also, 
\begin{eqnarray} \label{eq:h2b}\tilde\nabla_{X_i}h_*\xi_j&=&h_*\hat \nabla_{X_i}\xi_j+\alpha_h(f_*X_i, \xi_j)\nonumber\vspace{1ex}\\
&=&-F_*A^f_{\xi_j}X_i+h_*\nabla^\perp_{X_i}\xi_j+\frac{1}{r_1}\<\pi_1f_*X_i, \xi_j\>\nu_1+\frac{1}{r_2}\<\pi_2f_*X_i, \xi_j\>\nu_2
\nonumber\vspace{1ex}\\
&=&\delta_{ij}(-\lambda_i Z_i+\rho_4 h_*\xi_4+\frac{\lambda_i}{r_2}\nu_2),
\end{eqnarray}
where we have used (\ref{nxixii}). 

If $\kappa\lambda_i-k_1\neq 0$, that is,  $\<Z_i, Z_i\>\neq 0$,  define
$$
W_i = \tilde\nabla_{X_i}h_*\xi_i-\frac{\<\tilde\nabla_{X_i}h_*\xi_i,Z_i\>}{\<Z_i, Z_i\>}Z_i=\frac{-k_2\lambda_i}{\kappa \lambda_i-k_1}Z_i+\rho_4h_*\xi_4+\frac{\lambda_i}{r_2}\nu_2, \,\,\,1\leq i\leq 2.
$$
Then the vectors $F_*X_i, Z_i, h_*\xi_i$ and $W_i$ are pairwise orthogonal and the subspaces $V_i=\spa \{F_*X_i, Z_i, h_*\xi_i, W_i\}$, $1\leq i\leq 2$,  are  orthogonal to each other.
Using the second equations in  (\ref{rho}) and (\ref{ncd-imS-perp}), we obtain that
$\tilde \nabla_{X_i}W_j=0$
and
\be\label{eq:nablawi}\tilde \nabla_{X_i}W_i=\frac{k_1k_2(\lambda_i-\lambda_j)}{\kappa\lambda_i-k_1}h_*\xi_i, \,\,\,\,1\leq i\neq j\leq 2.\ee
It follows that the subspaces $V_1$ and $V_2$   are constant,
and that $\R_2^8$ also splits orthogonally as $\R_2^8=V_1\oplus V_2$. 

If $\kappa\lambda_i-k_1=0$, define
$$\zeta_i=\frac{\kappa}{2k_1k_2(\kappa\lambda_j-k_1)}(-2\kappa\tilde\nabla_{X_i}h_*\xi_i +(k_2-k_1)Z_i), \,\,\,1\leq i\neq j\leq 2.$$
Then $\<\zeta_i, \zeta_i\>=0$, $\<\zeta_i, Z_i\>=1$ and $\zeta_i\in \spa\{F_*X_i, h_*\xi_i\}^\perp$. Moreover, the subspaces $V_i=\spa \{F_*X_i, Z_i, h_*\xi_i, \zeta_i\}$, $1\leq i\leq 2$,  are  orthogonal to each other. Furthermore, since
\be\label{eq:nablazetai}
\tilde \nabla_{X_i}\zeta_i=\frac{k_1^2-k_2^2}{2k_1k_2}h_*\xi_i,
\ee
it follows that $V_1$ and $V_2$   are constant
and that $\R_2^8$ also splits orthogonally as $\R_2^8=V_1\oplus V_2$.

Since $\Gamma_{11}^2=\Gamma_{22}^1=0$, for each $x\in M^2$ there exists an isometry  $\psi\colon\, W=I_1\times I_2\to U_x$
of a product of open intervals $I_j\subset \R$, $1\leq j\leq 2$, onto a neighborhood of $x$,  such that $\psi_*\frac{\d}{\d s}=X_1$ 
and $\psi_*\frac{\d}{\d t}=X_2$, where $s$ and $t$ are the standard coordinates on $I_1$ and $I_2$, respectively.  Write $g=F\circ \psi$. In terms of the coordinates $(s,t)$, the fact that $\alpha_F(X_1, X_2)=0$
translates into 
$$\frac{\d^2 g}{\d s\d t}=0,$$
which implies that there exist smooth curves $\gamma_1\colon I_1\to V_1$ and $\gamma_2\colon I_2\to V_2$ such that $g=\gamma_1\times \gamma_2$.

If $\kappa\lambda_i-k_1\neq 0$, it follows from  (\ref{eq:zib}),  (\ref{eq:h1b}), (\ref{eq:h2b}) and (\ref{eq:nablawi}) that  $\gamma_i$ is a unit-speed space like curve  in $V_i$  with constant Frenet curvatures $\hat{k}^i_\ell$, $1\leq \ell\leq 3$, and Frenet frame $\{F_*X_i, \hat{Z}_i, h_*\hat{\xi}_i, \hat{W}_i\}$, where $\hat{Z}_i$, $\hat{\xi}_i$ and $\hat{W}_i$ denote the unit vectors in the direction of $Z_i$, $\xi_i$ and $W_i$, respectively.  Moreover, by (\ref{eq:zib}) and (\ref{eq:zinormb}) we have
$$(\hat{k}^i_1)^2=|\<Z_i,Z_i\>|=|(\lambda_i-\lambda_j)(\kappa\lambda_i-k_1)|,$$
whereas from (\ref{eq:h1}) and (\ref{eq:nablawi}) we obtain, respectively, that
$$(\hat{k}^i_2)^2=\frac{\kappa^2(\lambda_j-\lambda_i)^2\<\xi_i,\xi_i\>}{|\<Z_i, Z_i\>|}=\frac{\kappa^2|\lambda_j-\lambda_i|\lambda_i(1-\lambda_i)}{|\kappa\lambda_i-k_1|}$$
and
$$(\hat{k}^i_3)^2=\frac{k^2_1k^2_2(\lambda_i-\lambda_j)^2\<\xi_i, \xi_i\>}{(\kappa\lambda_i-k_1)^2|\<W_i, W_i\>|}=\frac{k_1k_2|\lambda_i-\lambda_j|}{|\kappa\lambda_i-k_1|},\,\,\,1\leq j\neq i\leq 2.$$

If $\kappa\lambda_i-k_1= 0$, it follows from  (\ref{eq:zib}),  (\ref{eq:h1b}), (\ref{eq:h2b}) and (\ref{eq:nablazetai}) that  $\gamma_i$ is a unit-speed space like curve  in $V_i$  with light-like curvature vector, constant Frenet curvatures $\tilde{k}^i_\ell$, $1\leq \ell\leq 2$, and Frenet frame $\{F_*X_i, {Z}_i, h_*\hat{\xi}_i, \zeta_i\}$, where  $\hat{\xi}_i$ is the unit vector in the direction of  $\xi_i$.  Moreover, from (\ref{eq:h1b}) we obtain that
$$(\tilde{k}^i_1)^2=\frac{k_1k_2(\kappa\lambda_j-k_1)^2}{\kappa^2},$$
whereas from  (\ref{eq:nablazetai}) it follows that
$$(\tilde{k}^i_2)^2=\<\xi_i, \xi_i\>\frac{(k_1^2-k_2^2)^2}{4k^2_1k^2_2}=\frac{(k_1-k_2)^2}{4k_1k_2}.
$$
Comparying with (\ref{eq:hatk1}), (\ref{eq:hatk2}) and (\ref{eq:hatk3}) in the first case, and with (\ref{eq:tildek1}) and (\ref{eq:tildek2}) in the second, we see that $\gamma_1$ and $\gamma_2$ are precisely, up to congruence, the curves given in Example~\ref{ex:2}.

Now observe that
$$\tilde{\pi}_2F_*\xi_i=h_*\pi_2 f_*\xi_i=h_*(f_*\xi_i+SX_i)=\lambda_iF_*X_i+h_*\xi_i,$$
whereas
$$\tilde{\pi}_2h_*\xi_i=h_*\pi_2\xi_i=h_*(f_*S^t\xi_i+T\xi_i)=(1-\lambda_i)(\lambda_iF_*X_i+h_*\xi_i),$$
where we have used that 
$$S^t\xi_i=S^tSX_i=R(I-R)X_i=\lambda_i(1-\lambda_i)X_i$$
and 
$$T\xi_i=TSX_i=S(I-R)X_i=(1-\lambda_i)SX_i=(1-\lambda_i)\xi_i.$$
On the other hand, 
$$\tilde{\pi}_2 Z_i=h_*\pi_2H+\frac{\lambda_i}{r_2}\pi_2\nu_2=\rho_4h_*\xi_4+\frac{\lambda_i}{r_2}\nu_2.$$
Since 
$$\tilde{\pi}_2 h_*\xi_4=h_*\pi_2\xi_4=h_*(f_*S^t\xi_4+T\xi_4)=h_*\xi_4,$$
we obtain that
$$\tilde{\pi}_2(\rho_4h_*\xi_4+\frac{\lambda_i}{r_2}\nu_2)=\rho_4h_*\xi_4+\frac{\lambda_i}{r_2}\nu_2.$$
If $\<Z_i, Z_i\>\neq 0$, it follows that  $\tilde{\pi}_2 W_i$ and $\tilde{\pi}_2 Z_i$ are colinear. Similarly, 
$\tilde{\pi}_2 \zeta_i$ and $\tilde{\pi}_2 Z_i$ are colinear if $\<Z_i, Z_i\>=0$. It follows that $\tilde {\pi}_2(V_i)$ is spanned by 
$$\lambda_iF_*X_i+h_*\xi_i\,\,\,\,\,\mbox{and}\,\,\,\,\,\rho_4h_*\xi_4+\frac{\lambda_i}{r_2}\nu_2.$$
Therefore, the subspaces  $\tilde {\pi}_2(V_1)$ and $\tilde {\pi}_2(V_2)$ (and hence also  $\tilde {\pi}_1(V_1)$ and $\tilde {\pi}_1(V_2)$)
are mutually orthogonal, thus the first (respectively, second) factor $\R_1^4$ in the decomposition $\R_1^4\times \R_1^4$ adapted to the 
product $\Hy_{k_1}^{3}\times \Hy_{k_2}^{3}$ splits orthogonally as $\R_1^4=\tilde{\pi}_1(V_1) \oplus \tilde{\pi}_1(V_2)$  (respectively, 
$\R_1^4=\tilde{\pi}_2(V_1) \oplus \tilde{\pi}_2(V_2)$). We conclude that  $g$ is (the restriction to $W$ of ) an isometric immersion as in Example \ref{ex:2}, and the conclusion follows as in the preceding case. 
\vspace{2ex}\qed

\section{The main result}

We are now in a position to state and prove our main result.

\begin{theorem}\po\label{thm:umb} Let $f\colon\,M^2\to \Q_{k_1}^{n_1}\times \Q_{k_2}^{n_2}$, $k_1+k_2\neq 0$,  be an umbilical non totally geodesic isometric immersion. 
Then one of the following possibilities holds:
\begin{itemize}
\item[$(i)$]  $f$ is an umbilical isometric immersion into a slice of $\Q_{k_1}^{n_1}\times \Q_{k_2}^{n_2}$;
\item[$(ii)$] there exist umbilical isometric immersions $f_i\colon\, M^2\to \Q_{\tilde k_i}^{n_i}$, $1\leq i\leq 2$, with $\tilde{k}_1=k_1\cos^2\theta$ and $\tilde k_2=k_2\sin^2\theta$ for some $\theta\in (0, \pi/2)$,  such that $f=(\cos \theta f_1,\sin\theta f_2)$;
\item[(iii)] after interchanging the factors, if necessary, we have $k_2=0$, $n_1\geq 3$, $n_2\geq 2$ and $f=j\circ \tilde f$, where $j\colon \Q_{k_1}^{3}\times \R^2\to \Q_{k_1}^{n_1}\times \R^{n_2}$ and $\tilde f\colon M^2\to \Q_{k_1}^{3}\times \R^2$ are isometric immersions such that $j$ is totally geodesic and $\tilde f(M^2)$ is an open subset of a surface as in Example \ref{ex:1};
\item[(iv)]  $k_i<0$ and  $n_i\geq 3$, $1\leq i\leq 2$,  and $f=j\circ \tilde f$, where 
$j\colon \Q_{k_1}^{3}\times \Q_{k_2}^{3}\to \Q_{k_1}^{n_1}\times \Q_{k_2}^{n_2}$ and  $\tilde f\colon M^2\to \Q_{k_1}^{3}\times \Q_{k_2}^{3}$ are isometric immersions such that $j$ is totally geodesic and $\tilde f(M^2)$ is an open subset of a surface as in Example \ref{ex:2};
\item[$(v)$] after possibly reordering the factors, we have $k_1>0$ (respectively, $k_1\leq 0$) and $f\circ \tilde\Pi=j\circ  \Pi\circ \tilde f$ (respectively, $f=j\circ  \Pi\circ \tilde f$), where $\tilde \Pi\colon\, \tilde M^2\to M^2$ is the universal covering of $M^2$, $\tilde f\colon\,\tilde M^2\to \R\times \Q_{k_2}^{2+\delta}$ (respectively, $\tilde f\colon\, M^2\to \R\times \Q_{k_2}^{2+\delta}$)  is  an umbilical isometric immersion with  $\delta\in \{0,1\}$, $j\colon\, \Q_{k_1}^{1}\times \Q_{k_2}^{2+\delta}\to \Q_{k_1}^{n_1}\times \Q_{k_2}^{n_2}$ is  totally geodesic  and $ \Pi\colon\,\R\times \Q^{2+\delta}_{k_2}\to  \Q_{k_1}^1\times \Q^{2+\delta}_{k_2}$ is  a locally isometric covering map (respectively, isometry).
\end{itemize}
\end{theorem}
\proof  
If $S$ vanishes everywhere on  $M^2$, then $f$ is as in  $(i)$  by Lemma $8.1$ in \cite{mt1}. If $\ker S=\{0\}$ at some point $x\in M^2$, then $f$ is as in $(ii)$, $(iii)$ or $(iv)$ by Lemma \ref{le:kerS=0}. Then, we are left with the case in which there is  an open subset ${\cal U} \subset M^2$ where $\ker S$ has rank one.  In this case, the argument in the proof of Theorem $1.4$ of \cite{mt1} applies and shows that $f$ is as in $(v)$. \qed

\vspace*{-1ex} {\renewcommand{\baselinestretch}{1}
\hspace*{-25ex}\begin{tabbing}
\indent  \=  Universidade Federal de S\~ao Carlos \\
\indent  \= Via Washington Luiz km 235 \\
\> 13565-905 -- S\~ao Carlos -- Brazil \\
\> e-mail: jlorjuelac@gmail.com \\
\hspace*{10ex} tojeiro@dm.ufscar.br
\end{tabbing}}


\begin{thebibliography}{99}

\bibitem{ckv} G. Calvaruso, D. Kowalczyk and J. Van der Veken, \ {\em On extrinsically symmetric hypersurfaces in $\Hy^n\times \R$.}   Bull. Aust. Math. Soc. \textbf{82} (2010),  390--400.

 \bibitem{ltv} J. H. Lira, R. Tojeiro and  F. Vit\'orio, \emph{ A Bonnet theorem for isometric immersions into products of space forms}, Archiv der Math. \textbf{95} (2010),  469--479.  
 
 \bibitem{mt1} B. Mendon\c ca and  R. Tojeiro, \ {\em Submanifolds of products of space forms},  Indiana Univ. Math. J. \textbf{62 (4)} (2013), 1283--1314.
 
 \bibitem{mt2} B. Mendon\c ca and  R. Tojeiro, \ {\em Umbilical submanifolds of $\Sf^n\times \R$},  Canadian J.  Math. \textbf{66 (2)} (2014), 400--428.
 
 \bibitem{mt3} B. Mendon\c ca and  R. Tojeiro, \ {\em Umbilical submanifolds of $\Hy^n\times \R$}, in preparation.
 
 \bibitem{al} Y. Nikolayevsky, \emph{Totally umbilical submanifolds of symmetric spaces}, Mat. Fiz. Anal. Geom. \textbf{1}
(1994), 314--357.

\bibitem{st} R. Souam and E. Toubiana,\ {\em Totally umbilic surfaces in homogeneous $3$-manifolds.}   Comment. Math. Helv. \textbf{84 (3)} (2009),  673--704.

\bibitem{vv} J. Van der Veken and L. Vrancken, \ {\em Parallel and semi-parallel hypersurfaces of $\Sf^n\times \R$.}   Bull. Braz. Math. Soc.  \textbf{39} (2008),  355--370.




\end{thebibliography}
\end{document}